\documentstyle{amsppt}
\magnification=1100
\NoBlackBoxes

\topmatter

\nopagenumbers
\font\titlefont=cmbx12\font\authorfont=cmr10
\baselineskip=13pt
\rightheadtext{SUBALGEBRAS AND DISCRIMINANTS OF ANTICOMMUTATIVE ALGEBRAS}
\leftheadtext{E.A. TEVELEV}
\topmatter
\centerline{\titlefont
SUBALGEBRAS AND DISCRIMINANTS}
\smallskip
\centerline{\titlefont
OF ANTICOMMUTATIVE ALGEBRAS}
\vskip 1.2pc
\centerline{\authorfont
E.A. TEVELEV$^{\ast}$}
\vskip 1.1pc
\centerline{Moscow State University}
\centerline{Vorobievy Gory, Moscow, 119899, Russia}
\smallskip
\centerline{tevelev\@ium.ips.ras.ru}

\thanks{$^\ast$ This research was carried out with 
the financial support of the CRDF Foundation
(grant RM1-206).}\endthanks

\abstract
The paper deals with the configuration of subalgebras in generic
$n$-dimensional  $k$-argument anticommutative algebras and
``regular'' anticommutative algebras.
\endabstract

\endtopmatter
\def\cal#1{{\fam2 #1}}
\def\C{{\Bbb C}}
\def\A{{\cal A}}
\def\D{{\cal D}}
\def\E{{\cal E}}
\def\L{{\cal L}}
\def\SS{{\cal S}}
\def\gl{\goth{gl}}

\def\Gr{\text{Gr}}
\def\X{\text{{\bf X}}}
\def\GL{\text{GL}}
\def\SL{\text{SL}}
\def\codim{\text{codim}}
\def\Tr{\mathop{\text{Tr}}}

\def\Pic{\text{Pic}}
\def\dualA{{\A_{n,k}^0}^{\!\!\!\!*}}
\def\eps{\varepsilon}
\def\a{\alpha}
\def\schoose#1#2{\left[\matrix #1\cr #2\endmatrix\right]}
\def\fchoose#1#2{\left\{\matrix #1\cr #2\endmatrix\right\}}
\def\ds{\displaystyle}

\document
\head Introduction\endhead

Let $V=\C^n$. 
We fix an integer $k$, $1<k<n-1$. 
Let $\A_{n,k}=\Lambda^kV^*\otimes V$ be the vector space
of $k$-linear anticommutative maps from
$V$ to $V$. We identify the points of $\A_{n,k}$ 
with the corresponding algebras, that is, we assume that
$A\in \A_{n,k}$ is the space $V$ 
equipped with the structure of $k$-argument
anticommutative algebra.

\subhead 1. Subalgebras in generic algebras\endsubhead
Subalgebras in generic algebras with $k=2$ 
were studied in \cite{11}.
The following theorem is a generalization of these results.

\proclaim{Theorem 1}
Let $A\in\A_{n,k}$ be a generic algebra. Then\par
{\rm(i)} every $m$-dimensional subspace is a subalgebra if $m<k$,\par
{\rm(ii)} $A$ contains no $m$-dimensional subalgebras
with $k+1<m<n$,\par
{\rm(iii)} the set of $k$-dimensional subalgebras
is a smooth irreducible $(k-1)(n-k)$-dimensional
subvariety in the Grassmanian $\Gr(k,A)$,\par
{\rm(iv)} there are finitely many $(k+1)$-dimensional subalgebras, and their
number is
$$\sum_{
\matrix
\scriptstyle n-k-1\ge\mu_1\ge\ldots\ge\mu_{k+1}\ge 0\cr
\scriptstyle n-k-1\ge\lambda_1\ge\ldots\ge\lambda_{k+1}\ge 0\cr
\scriptstyle\mu_1\le\lambda_1,\ldots,\mu_{k+1}\le\lambda_{k+1}
\endmatrix
}\!\!\!\!\!\!\!\!\!
(-1)^{|\mu|}
\frac{(\lambda_1+k)!(\lambda_2+k-1)!\ldots\lambda_{k+1}!}
{(\mu_1+k)!(\mu_2+k-1)!\ldots\mu_{k+1}!}
(|\lambda|-|\mu|)!
\left|\frac{1}{(i-j+\lambda_j-\mu_i)!}\right|^2_{i,j=1,\ldots,k+1},$$
where $|\lambda|=\lambda_1+\ldots+\lambda_{k+1}$, 
$|\mu|=\mu_1+\ldots+\mu_{k+1}$, $1/N!=0$ if $N<0$,\par
{\rm(v)} $A$ contains a $(k+1)$-dimensional subalgebra,\par
{\rm(vi)} if $k=n-2$, then the number of $(k+1)$-dimensional subalgebras
is equal to $$2^n-(-1)^n\over 3.$$
\endproclaim

This theorem will be proved in \S1. Here we give a scetch of the proof.
We begin with the following general situation.

Assume that $G$ is a connected reductive group, 
$T$  is a fixed maximal torus,
and $B$ is a fixed Borel subgroup, $T\subset B\subset G$,
$B_-$ is the opposite Borel subgroup,
$P$ is a parabolic subgroup, $P\supset B_-$,
$\X(T)$ is the lattice of characters of $T$,
and $\lambda\in\X(T)$  is the dominant weight.
Consider the vector bundle $\L_\lambda=G\times_PU_\lambda$ over $G/P$, where
$U_\lambda$ is the irreducible $P$-module with highest weight $\lambda$.
By the Borel--Weil--Bott theorem (see~\cite{7}), 
$V_\lambda=H^0(G/P, \L_\lambda)$ is an irreducible $G$-module with highest
weight $\lambda$. 

\proclaim{Lemma 1} Let $s\in V_\lambda$ be a generic global section. Then\par
{\rm(i)} if $\dim U_\lambda>\dim G/P$, the scheme of zeros $Z_s$ is empty,\par
{\rm(ii)} if $\dim U_\lambda\le\dim G/P$, either $Z_s$ is empty
or $s$ intersects the zero section $\L_\lambda$ transversally and $Z_s$ 
is a smooth unmixed subvariety of codimension
$\dim U_\lambda$,\par
{\rm(iii)} if $\dim U_\lambda=\dim G/P$, the geometric number of points
in $Z_s$ is equal to the highest Chern class of $\L_\lambda$.\par
\endproclaim

The proof will be given in \S1. 
Here we show how this lemma can be used to prove the theorem.
The $\GL_n$-module $\A_{n,k}$ is a sum of two irreducible submodules:
$$\A_{n,k}=\A_{n,k}^0\oplus\tilde\A_{n,k}.\eqno(1)$$
Here $\tilde\A_{n,k}$ is isomorphic to $\Lambda^{k-1}V^*$:
we assign to every $(k-1)$-form
$\omega$ the algebra with multiplication
$$[v_1,\ldots,v_k]=\sum_{i=1}^k(-1)^{i-1}
\omega(v_1,\ldots,\hat v_i,\ldots,v_k)v_i.$$

Note tha every subspace of this algebra is a subalgebra.
Hence, the lattice of subalgebras of $A\in\A_{n,k}$
coincides with that of $A^0$ (the zero component of $A$), 
where $A\mapsto A^0$ is the $\GL_n$-equivariant
projector on the first summand in (1).
Algebras in $\A_{n,k}^0$ will be called {\it zero trace algebras\/},
since $A\in\A_{n,k}^0$ if and only if the $(k-1)$-form
$\Tr[v_1,\ldots,v_{k-1},\cdot]$ is equal to zero. 
Hence, the theorem will be proved
once we have proved
it for generic algebras in $\A_{n,k}^0$.

We choose a basis $\{e_1,\ldots,e_n\}$ in $V$, identify
$\GL_n$ with the group of matrices, 
consider the standard diagonal maximal torus $T$, 
and take for $B$ and $B_-$ the subgroups of upper- and lower-triangular
matrices. We fix an  $m\ge k$. 
Consider the parabolic subgroup of matrices
$$P=\left(\matrix
A&0\cr
*&B\cr\endmatrix\right),\eqno(2)$$
where $B$ is an $m\times m$-matrix and $A$ is an $(n-m)\times(n-m)$-matrix.
Then $G/P$ coincides with $\Gr(m,V)$. 
Consider the vector bundle  $\L=\Lambda^k\SS^*\otimes V/\SS$
on $G/P$, where $\SS$ is the tautological bundle and 
$V/\SS$ is the factor-tautological bundle. 
The assumptions of the lemma are fulfilled, since
$\L=\L_\lambda$, where $\lambda$ is the highest weight of $\A_{n,k}^0$. 
Therefore, 
$\A_{n,k}^0=H^0(\Gr(m,V), \Lambda^k\SS^*\otimes V/\SS)$.
Let $A\in\A_{n,k}^0$, and let $s_A$ be the corresponding
global section.
Then $(Z_{s_A})_{red}$ coinsides with the variety
of $m$-dimensional subalgebras of $A$.

Let us return to the theorem.
Assertion (i) is obvious.
Assertion (ii) follows from assertion (i) of the lemma.
Assertion (ii) of the lemma implies that if every $k$-argument
anticommutative algebra $A$ contains a $k$-dimensional subalgebra,
then the variety of $k$-dimensional subalgebras of a generic algebra
is a smooth unmixed $(k-1)(n-k)$-dimensional subvariety in $\Gr(k,A)$.
We claim that any $(k-1)$-dimensional subspace $U$
can be included in a $k$-dimensional subalgebra. 
Themultiplication in the algebra defines a linear map from $V/U$ to $V/U$. 
Let $v+U$ be a non-zero eigenvector.
It is obvious that $\C v\oplus U$ is a $k$-dimensional subalgebra. 
To prove the irreducibility of the variety of $k$-dimensional
subalgebras we use the Koszul complex. Assertions (iv) and (vi) of the theorem
follow from assertion (iii) of the lemma and explicit calculations
in the Chow ring of $\Gr(k+1,V)$. 
Assertion (v) of the theorem requires additional calculations. 

It should be noted that Lemma~1 cannot be strengthened to the point
where the non-emptyness and the irreducibility of the scheme of
zeros in Theorem 1 could be established apriori, as the following
example shows.
Consider the vector bundle $S^2\SS^*$ on $\Gr(k, 2n)$.
The dimension of a fibre does not exceed the dimension of
the Grassmanian as $k\le {\ds 4n-1\over\ds 3}$, but a generic section
(that is, a non-degenerate quadratic form in $\C^{2n}$) 
has a zero (that is, a $k$-dimensional isotropic subspace)
only if $k\le n$. For $k=n$ the scheme of zeros is a reducible variety of dimension
$\ds n(n-1)\over \ds2$ with two irreducible components that
correspond to two families of maximal isotropic subspaces on an even-dimensional quadric.

An essential drawback of the theorem is the fact that
it does not enable us to study the structure of subalgebras of any particular algebra.
The purpose of the remaining part of the paper is to correct this situation.

\subhead 2. $D$-regular algebras\endsubhead
Let $\dualA$ be the $\GL_n$-module dual to $\A_{n,k}^0$,
and let $S_D$ be the closure of the orbit of the highest
vector, $S_D\subset\dualA$. Let $PS_D\subset P\dualA$
be its projectivization, let 
$P\D\subset P\A_{n,k}^0$ be the subvariety projectively dual
(see~\cite{2}) to the subvariety $PS_D$, 
and let $\D\subset \A_{n,k}^0$ be the cone over it.
Then $\D$ is called the {\it $D$-discriminant subvariety}.
The algebras $A\in\D$ are said to be {\it $D$-singular}.
The algebras $A\not\in\D$ are said to be {\it $D$-regular}.

\proclaim{Theorem 2}
{\rm(i)} $\D$ is a hypersurface.\par
{\rm(ii)} Let $A$ be a $D$-regular algebra. Then the set of 
$k$-dimensional subalgebras of $A$ is a smooth irreducible 
$(k-1)(n-k)$-dimensional subvariety in $\Gr(k,A)$.\par
{\rm(3)} Let $k=n-2$. Then the degree of $\D$ is equal to
$${(3n^2-5n)2^n-4n(-1)^n}\over{18}.\eqno(3)$$
\endproclaim

Hence, the $D$-singularity of $A$ is determined
by the vanishing of the $SL_n$-invariant polynomial $D$
that defines~$\D$. 
This polynomial is called the $D$-discriminant.
Theorem~2 will be proved in \S2. 
Here we give a sketch of the proof.
The fact that $\D$ is a hypersurface follows immediately
from the results of \cite{9}.
Assertion (ii) can be deduced from the corresponding assertion
of Theorem~1 by an easy calculation with differentials.
Assertion (iii) requires some comments.
It is easy to show that $PS_D$ coincides with the variety of incomplete flags
$0\subset V_1\subset V_2\subset \Bbb C^n$, 
where $\dim V_1=1$ and $\dim V_2=n-k$, in the ``Pl\"ucker'' embedding.
Hence,one should not hope
to find a closed formula for the degree of
$\D$, since there is no good formula even for the degree of the variety projectively dual to the 
Grassmannian
in the Pl\"ucker embedding.
(See~\cite{10} and the formula for the variety
projectively dual to the variety of complete flags in~\cite{8}).
The following theorem is a generalization of formula (3).

\proclaim{Теорема 2'}
Let $V^*$ be an irreducible $\SL_n$-module with highest
weight $(a-1)\varphi_1+\varphi_2$. Then the variety
$\D\subset PV$ projectively dual to the projectivization
of the orbit of the highest vector is a hypersurface of degree
$${(n^2-n)a^{n+1}-(n^2+n)a^{n-1}-2n(-1)^n\over (a+1)^2}.$$
\endproclaim

\subhead 3. $E$-regular algebras\endsubhead 
We define the $E$-discriminant and $E$-regularity
only for $(n-2)$-argument $n$-dimensional anticommutative algebras. 
Let $\A=\A_{n,n-2}^0$.
Consider the projection $\pi:\,\Gr(n-1,V)\times P\A\to P\A$
on the second summand and the incidence subvariety
$Z\subset\Gr(n-1,V)\times P\A$ that consists of $S\subset PA$,
where $S$ is a subalgebra in $A$.
Let $\tilde\pi=\pi|_Z$.
By Theorem~1, we have $\tilde\pi(Z)=\A$. 
Let $\tilde\E\subset Z$ be the set of criticalpoints of $\tilde\pi$,
let $P\E=\tilde\pi(\tilde\E)$ be the set of critical values of
$\tilde\pi$, and let $\E\subset\A$ be the cone over $P\E$.
Then $\E$ is called the {\it $E$-discriminant subvariety}.
The algebras $A\in\E$ are said to be {\it $E$-singular}.
The algebras $A\not\in\E$ are said to be {\it $E$-regular}.

\proclaim{Theorem 3}
{\rm(i)} $\E$ is an irreducible hypersurface.\par
{\rm(ii)} Let $A$ be an $E$-regular algebra. 
Then $A$ has precisely
$$2^n-(-1)^n\over 3$$ 
$(n-1)$-dimensional subalgebras.\par
{\rm(iii)} The map $\tilde\pi:\,\tilde\E\to P\E$ is birational.
\endproclaim

Hence, the $E$-singularity of $A$ is determined by the vanishing of the 
$SL_n$-invariant polynomial that defines~$\E$. 
This polynomial is called the {\it $E$-discriminant}. 
Assertion (iii) can be formulated as follows: a generic
$E$-singular algebra has precisely one ``critical'' 
$(n-1)$-dimensional subalgebra.

Theorem 3 will be proved in \S3. 
Here we give a scetch of the proof.
We deduce assertion (i) from assertion (iii) and the irreducibility of 
$\tilde\E$ by calculating the dimension of $\tilde\E$.
Assertion (ii) can be proved by simplecalculation with differentials.
Hence, we have only to prove assertion (iii).

\subhead 4. Regular $4$-dimensional anticommutative algebras\endsubhead
An $(n-2)$-argument $n$-dimensional anticommutative algebra is said to be 
{\it regular} if it is $D$-regular and $E$-regular.
In this subsection we consider $2$-argument
$4$-dimensional algebras.
The corresponding generic algebras
were studied in \cite{11}. П
Here we formulate the statements
on $4$-dimensional generic algebras that are valid for all regular algebras.

\proclaim{Theorem 4}
Let $A$ be a $4$-dimensionalregular anticommutative algebra.
Then the following assertions are valid.\par
{\rm(i)} $A$ has precisely five $3$-dimensional subalgebras.
The set of these subalgebras
is a generic configuration of five hyperplanes.
In particular, $A$ has a pentahedral normal form,
that is, it can be reduced by a transformation
that belongs to  $\GL_4$ to an algebra such that the set of its five subalgebras is
a Sylvester pentahedron
$x_1=0$, $x_2=0$, $x_3=0$, 
$x_4=0$, $x_1+x_2+x_3+x_4=0$.\par
{\rm(ii)} $A$ has neither one- nor two-dimensional ideals.\par
{\rm(iii)} The set of two-dimensional
subalgebras of $A$ is a del Pezzo surface of degree~$5$ 
(a blowing up of $\Bbb P^2$ at four generic points.\par
{\rm(iv)} $A$ has precisely $10$ fans, that is,
flags $V_1\subset V_3$ of
$1$-dimensional and $3$-dimensional
subspaces such that every intermediate subspace
$U$, $V_1\subset U\subset V_3$,
is a two-dimensional subalgebra.
\endproclaim

The proof will be given in \S4. 
There are examples showing that the other statements on generic algebras
proved in  \cite{11}
are not, generally speaking,
valid for regular algebras.
There are examples of regular algebras
that have more than two
commutative subalgebras,
that have three-dimensional ideals,
whose associated cubic hypersurface is not smooth (not irreducible), and so on.

The author is grateful to E.B.~Vinberg
for useful discussions
and the simplification of some proofs.

\head \S1. Subalgebras in generic algebras\endhead

\subhead Proof of Lemma 1\endsubhead
Assertion (i) is obvious.
It was proved, for example, in \cite{11}, Lemma 1.1.

We now prove assertion (ii). Let $\dim U_\lambda\le\dim G/P$,
and assume that every global section has a zero.
We have to prove that a generic global section
$s$ intersects the zero section of $\L_\lambda$ 
transversally. This will imply, in particular, that 
$Z_s$ is a smooth unmixed subvariety of codimension $\dim U_\lambda$. 
For simplicity we suppress the index~$\lambda$. 
Consider $G/P\times V$ and $Z\subset G/P\times V$, 
$Z=\{(x,s)\,|\,x\in (Z_s)_{red}\}$.
Since $G/P$ is homogeneous and $Z$ is invariant,
it follows that $Z$  is obtained by spreading
the fibre $Z_e=\{s\in V\,|\,s(eP)=0\}$ by the group $G$.
Since $U$ is irreducible, we have $\dim Z_e=\dim V-\dim U$.
Hence, $Z$ is a smooth irreducible subvariety
of dimension $\dim V+\dim G/P-\dim U$.

Let $\pi:\,Z\to V$ be the restriction to $Z$ 
of the projection of $G/P\times V$ on the second summand.
By assumption, $\pi$ is a surjection.
By Sard's lemma for algebraic varieties (see~\cite{4}),
for a generic point $s\in V$ and any point $(x,s)$ 
in $\pi^{-1}s$ the differential $d\pi_{(x,s)}$ is surjective. 
We claim that $s$ has transversal intersection with the zero section.
Indeed, $\L=(G/P\times V)/Z$
is regarded as a vector bundle over $G/P$. 
The zero section of $\L$ is identified with $Z/Z$.
The section $s$ is identified with $(G/P\times\{s\})/Z$.
Hence, it is sufficient to prove that $G/P\times\{s\}$ 
is transversal to $Z$,
which is equivalent to the following assertion:
$d\pi_{(x,s)}$ is surjective for all $(x,s)\in Z$.
Now assertion (iii) of the lemma follows from the standard intersection 
theory (see~\cite{5}): if the scheme of zeros of a generic global
sectionis empty, then
$\L$ contains a trivial one-dimensional
subbundle, and its highest Chern class is zero.
If it is non-empty,
we use assertion (ii).\qed

\subhead Proof of Theorem 1\endsubhead
Assertion (i) is obvious.
To prove the other assertions we use
the realization of
$\A_{n,k}^0$ as $H^0(\Gr(m,V), \Lambda^k\SS^*\otimes V/\SS)$.
Assertions (ii) and (iii) follow immediately
from Lemma~1. We have only to prove that the variety
of $k$-dimensional subalgebras is irreducible.
Assume that a section $s$ of the bundle $\L=\Lambda^k\SS^*\otimes V/\SS$ 
over $\Gr(k,V)$corresponding to $A$ has transversal intersection
with the zero section. Then the Koszul complex
$$0\rightarrow\Lambda^{n-k}{\L}^*\buildrel s\over \rightarrow\ldots
\buildrel s\over \rightarrow\Lambda^2{\L}^*\buildrel s\over \rightarrow
{\L}^*\buildrel s\over \rightarrow{\cal O}
\rightarrow{\cal O}_{Z(s)}\rightarrow0$$
is exact (see~\cite{5}).

Note that $\Lambda^p{\L}^*$ is isomorphic to the bundle
$S^p\Lambda^k\SS\otimes \Lambda^p(V/\SS)^*$. 
This is a homogeneous bundle over $G/P$ 
of the form $\L_\mu$ (see Introduction), where
$\mu=-\eps_1-\ldots-\eps_p+\eps_{n-k+1}+\ldots+\eps_n$ and 
$\eps_i$ are the weights of the diagonal torus
in the tautological representation.
Note that $\mu+\rho$
(where $\rho$ is the half-sum of the positive roots)
is singular (belongs to the wall of the Weil chamber)
for any $p$, $1\le p\le n-k$.
By the Borel--Weil--Bott theorem,
$H^*(\Gr(k,V), \Lambda^p{\L}^*)=0$ for $1\le p\le n-k$.
Hence, $H^0(Z(s), {\cal O}_{Z(s)})=H^0(\Gr(k, V), {\cal O})=\C$,
as was to be shown.

We postpone the proof of assertion (v) till the end of this section
and consider assertions (iv) and (vi), that is,
we shall calculate the highest Chern class of the bundle
$\Lambda^k\SS^*\otimes V/\SS$ over $\Gr(k+1, V)$. 

We use the standard notation, facts, and formulae from the Schubert
calculus (cf.~\cite{5}. 
The letters $\lambda$ and $\mu$ allways denote Young
diagrams
in the rectangle with $k+1$ rows
and $n-k-1$ columns.
It is well known that such diagrams
parametrize the basis of the Chow ring of
$\Gr(k+1, V)$. The cycle corresponding to $\lambda$
is denoted by $\sigma_\lambda$, 
$\sigma_\lambda\in A^{|\lambda|}(\Gr(k+1, V))$. (We grade the Cjow ring
by the codimension, $|\lambda|=\lambda_1+\ldots+\lambda_{k+1}$, where
$\lambda_i$ is the length of the $i$th row of $\lambda$.)
We need the total Chern class of the bundles $\SS$ and $V/\SS$:
$$c(\SS)=1-\sigma_1+\sigma_{1,1}-\ldots+(-1)^{k+1}\sigma_{1,\ldots,1},
\qquad
c(V/\SS)=1+\sigma_1+\sigma_2+\ldots+\sigma_{n-k-1}.
$$

We begin by calculating the total Chern class of $\SS\otimes V/\SS$.
The standard formula for the total Chern class of the tensor products of two bundles implies that
$$c(\SS\otimes V/\SS)=\sum_{\mu\subset\lambda}d_{\lambda\mu}
\Delta_{\tilde\mu}(c(\SS))\Delta_{\lambda'}(c(V/\SS)),\eqno(4)$$
where
$$d_{\lambda\mu}=\left|{\lambda_i+k+1-i
\choose \mu_j+k+1-j}\right|_{1\le i,j\le k+1},\quad
\Delta_\lambda(c)=|c_{\lambda_i+j-i}|,$$ 
$\lambda'=(n-k-1-\lambda_{k+1},n-k-1-\lambda_k,\ldots,n-k-1-\lambda_1)$,
and $\tilde\mu$ is the diagram obtained from $\mu$ by transposition.
We shall use the fact that
$\Delta_{\tilde\mu}(c(\SS))=\Delta_{\mu}(s(\SS))$,
where $s(E)$ is the Segre class of $E$.
Another fact is that 
$\Delta_{\lambda}(c(V/\SS))=\sigma_{\lambda}$,
$\Delta_{\lambda}(s(\SS))=(-1)^{|\lambda|}\sigma_{\lambda}$
(since $s(\SS)=1-\sigma_1+\sigma_2+\ldots+(-1)^{n-k-1}\sigma_{n-k-1}$).
Therefore, (4) can be written as
$$c(\SS\otimes V/\SS)=\sum_{\mu\subset\lambda}d_{\lambda\mu}(-1)^{|\mu|}
\sigma_\mu\sigma_{\lambda'}.\eqno(5)$$

To calculate the highest Chern class of
$\L=\Lambda^k\SS^*\otimes V/\SS$,we use the formula
$\L=\Lambda^{k+1}\SS^*\otimes(\SS\otimes V/\SS)$. 
The total Chern class of the first factor is equal to
$c(\Lambda^{k+1}\SS^*)=1+\sigma_1$. The total Chern class of the second is given by (5). 
Hence, the highest Chern class of
$\L$ is
$$c_{top}(\L)=\sum_{\mu\subset\lambda}d_{\lambda\mu}(-1)^{|\mu|}
\sigma_\mu\sigma_{\lambda'}\sigma_1^{|\lambda|-|\mu|}.\eqno(6)$$

The last result of the Schubert calculus that we need is the exact formula for the
degree of a product of two cycles.
In our case this can be written as
$$\sigma_\mu\sigma_{\lambda'}\sigma_1^{|\lambda|-|\mu|}=
\deg(\sigma_\mu\sigma_{\lambda'})=(|\lambda|-|\mu|)!
\left|\frac{1}{(i-j+\lambda_j-\mu_i)!}\right|_{i,j=1,\ldots,n},\eqno(7)$$
where $1/N!=0$ if $N<0$.
The formula in assertion (iv) of Theorem~1
can be obtained from (6) and (7) by a slight modification
of the determinant
in the formula for $d_{\lambda\mu}$.

It remains to verify the formula in the assertion (vi)
of Theorem~1. Let
$k=n-2$. Then the formula in assertion (iv) can be written as
$$\sum_{0\le i\le j\le k+1}
(-1)^i\frac{(1+k)!(1+(k-1))!\ldots(1+(k-j+1))!(k-j)!\ldots0!}
{(1+k)!(1+(k-1))!\ldots(1+(k-i+1))!(k-i)!\ldots0!}
(j-i)!\,{\det}^2A,\eqno(8)$$
where
$$A=\left(\matrix
X&0&0\cr
*&Y&0\cr
*&*&Z\cr
\endmatrix\right),$$
$X$ is an $i\times i$-matrix, 
$Y$ is a $(j-i)\times (j-i)$-matrix, and
$Z$ is a $(k+1-j)\times (k+1-j)$-matrix.
$X$ and $Z$ are lower-triangular matrices with $1$s on the diagonal and $Y$ is given by
$$Y=\left(\matrix
1   &1   &0   &0&\cdots&0\cr
1/2!&1   &1   &0&\cdots&0\cr
1/3!&1/2!&1   &1&\cdots&0\cr
1/4!&1/3!&1/2!&1&\cdots&0\cr
\vdots&\vdots&\vdots&\vdots&\ddots&\vdots\cr
1/(j-i)!&1/(j-i-1)!&1/(j-i-2)!&1/(j-3)!&\cdots&1\cr
\endmatrix\right).$$
It is easy to verify that $\det Y=1/(j-i)!$,
which enables us to rewrite (8) as
$$
\aligned
&\sum_{i\le j}
(-1)^i\frac
{(1+(k-i))!\ldots(1+(k-j+1))!}
{(k-i)!\ldots(k-j+1)!}
/(j-i)!=
\sum_{i\le j}
(-1)^i{k+1-i\choose j-i}=\cr
&\sum_i(-1)^i2^{k+1-i}={2^{k+2}-(-1)^{k+2}\over 3}=
{2^n-(-1)^n\over 3},\cr
\endaligned$$
as was to be shown.

It remains to prove assertion (v) of Theorem~1.
We have to prove that every
$A\in\A_{n,k}$ has a $(k+1)$-dimensional subalgebra.
We fix a basis $\{e_1,\ldots,e_n\}$ in $V$
and consider the subspace
$U=\langle e_{n-k}, \ldots, e_n\rangle$. 
Let $M\subset\A_{n,k}$ be the subspace that consists of the algebras for which
$U$ is a $(k+1)$-dimensional subalgebra.
We claim that $\A_{n,k}=\GL_n\cdot M$. 
It is sufficient to prove that the differential of the canonical morphism
$\phi:\,\GL_n\times M\to\A_{n,k}$ 
is surjective at a point $(e, A)$. 
Consider the algebra $A\in M$ in which
$[e_{n-k},\ldots,\hat e_i,\ldots, e_n]=e_i$
for all $n-k\le i\le n$ and the other products are zero.
We claim that $d\phi$ is surjective at $(e, A)$. 
Consider the map $\pi:\,\A_{n,k}\to\A_{n,k}/M $.
It is sufficient to verify that
$\pi\circ d\phi_{(e, A)}(\gl_n, 0)=\A_{n,k}/M\simeq\Lambda^kU^*\otimes V/U$.
But this is obvious, since multiplication in
$d\phi_{(e, A)}(E_{ji}, 0)$  for
$n-k\le i\le n$, $1\le j\le n-k-1$
($E_{ji}$ is the matrix identity)
is given by
$[e_{n-k},\ldots,\hat e_i,\ldots, e_n]=e_j$
with the other products equal to zero.
This completes the proof of the theorem.\qed

\head \S2. $D$-regular algebras\endhead

\subhead Proofs of Theorems 2 and 2'\endsubhead
The irreducible representations of semisimple groups
for which the variety dual to the projectivization
of the orbit of the highest vector
is not a hypersurface were found in~\cite{9}. 
Since this list does not contain  $\dualA$
(as an $\SL_n$-module), $\D$ is a hypersurface.
Hence, the discriminant $D$ is well defined.

We prove (ii). Let $A$ be a $D$-regular algebra.
We use the arguments
in the proof of assertion (ii) of Lemma~1 and assertion (iii) of Theorem~1.
According to these calculations, it is sufficient
to verify that if $U$ is a $k$-dimensional subalgebra of $A$,
then the map
$$\psi:\,\gl_n\to\Lambda^kU\otimes A/U,\quad
\psi(g)(v_1\wedge\ldots\wedge v_k)=
g[v_1,\ldots,v_k]-[gv_1,\ldots,v_k]-\ldots-[v_1,v_2,\ldots,gv_k]+U$$
is surjective.
Assume the contrary.
Then there is a hyperplane $H\supset U$ 
such that the image of $\psi$
lies in $\Lambda^kU\otimes H/U$.
Consider a non-zero algebra $\tilde A$ in $S_D\subset\dualA$
such that
$[U^\perp,V^*,\ldots,V^*]=0$, $[V^*,\ldots,V^*]\subset H^\perp$,
where $U^\perp$ and $H^\perp$ are the annihilators of $U$ and $H$ in
$\dualA$. (Such an algebra is unique up to a scalar.)
Then $\tilde A$ annihilates $[\gl_n, A]$, which is equivalent
to the fact that $A$ annihilates $[\gl_n, \tilde A]$, that is,
the tangent space to $S_D$ at $\tilde A$. 
This means that $A$ lies in $\D$, that is, it is a $D$-singular algebra.

In what follows we assume that $k=n-2$. 
Since assertion (iii) of Theorem~2
follows from Theorem~2' (with $a=2$),
it is sufficient to prove Theorem~~2'.
We have to calculate the degree of~$\D$.
We use Kleiman's formula (see~\cite{2})
for the degree of the dual variety:
if $Z$ is a smooth projective $l$-dimensional variety in $\Bbb P^{n-1}$ 
and ${\cal L}={\cal O}_{\Bbb P^{n-1}}(1)|_Z$, then
$$\deg(\check Z)=\sum_{i=0}^l(i+1)\int_Zc_{l-i}(\Omega_Z^1)c_1({\cal L})^i,$$
where $\check Z$ is the projectively dual variety.

In the present case $Z=G/P$, where $G=\GL_n$ and
$P\subset G$ is the parabolic subgroup of matrices
$$\left(\matrix
*&0&0&\ldots&0\cr
*&*&0&\ldots&0\cr
*&*&*&\ldots&*\cr
\vdots&\vdots&\vdots&\ddots&\vdots\cr
*&*&*&\ldots&*\cr
\endmatrix\right).\eqno(9)$$
Assume that $T\subset G$ is the diagonal torus,
$B$ is the Borel subgroup of lower-triangular matrices,
$x_1,\ldots,x_n$ are the weights of the tautological representation,
$X(T)$ is the lattice of characters of $T$, 
$S$ is the symmetric algebra of $X(T)$ (over~$\Bbb Q$),
$W\simeq S_n$ is the Weil group of $G$,
and $W_P\simeq S_{n-2}$ is the Weil group of~$P$.
It is well known (see~\cite{1}) that the map $c:\,X(T)\to\Pic(G/B)$
that assigns to~$\lambda$ the first Chern class of the invertible sheaf
${\cal L}_\lambda$ (see the Introduction)can be extended
to a surjective homomorphism
$c:\,S\to A^*(G/B)$ in the (rational) Chow ring,
and its kernel coincides with $S^W_+S$.
The projection $\alpha:\,G/B\to G/P$ induces an embedding 
$\alpha^*:\,A^*(G/P)\to A^*(G/B)$.
Theimage coincides with the subalgebra
of $W_P$-invariants.
Hence, $A^*(G/P)=S^{W_P}/S^{W}_+S^{W_P}$. 
We denote the homomorphism
$S^{W_P}\to A^*(G/P)$ by the same letter~$c$.

To apply Kleiman's formula we need
$c_1(\cal L)$ (which is equal to $c(ax_1+x_2)$) and the total Chern class of
$\Omega_Z^1$, which is equal to
$$c(\Omega_Z^1)=c\left((1-x_1+x_2)\prod_{i=3,\ldots,n}(1-x_1+x_i)
\prod_{i=3,\ldots,n}(1-x_2+x_i)\right).$$
(This can be shown by standard arguments using the filtration of~$\Omega_Z^1$.
See also~\cite{8}).

Let $\a_1,\ldots,\a_{n-2}$ be the elementary symmetric polynomials in
$x_3,\ldots,x_{n}$. Then\break
$S^{W_P}=\Bbb Q[x_1,x_2,\a_1,\ldots,\a_{n-2}]$,
and the ideal $S^W_+S^{W_P}$ is generated by
$$\matrix
&x_1+x_2+\a_1,\quad x_1x_2+x_1\a_1+x_2\a_1+\a_2,\cr
&x_1x_2\a_1+x_1\a_2+x_2\a_2+\a_3,\ \ldots,\ 
x_1x_2\a_{n-4}+x_1\a_{n-3}+x_2\a_{n-3}+\a_{n-2},\cr
&x_1x_2\a_{n-3}+x_1\a_{n-2}+x_2\a_{n-2},\quad x_1x_2\a_{n-2}.
\endmatrix$$
Hence,
$\a_i=(-1)^i(x_1^i+x_1^{i-1}x_2+\ldots+x_2^i) \mod S^W_+S^{W_P}$, and
$A^*(G/P)$  is isomorphic to the quotient ring
$\Bbb Q[x_1, x_2]/\langle f_1, f_2\rangle$, where
$f_1=x_1^{n-1}+x_1^{n-2}x_2+\ldots+x_2^{n-1}$ and
$f_2=x_1^n$. Note that $f_1,f_2$ is the Gr\"obner basis of the ideal
$\langle f_1, f_2\rangle$ with respect to the ordering
$x_2>x_1$ (see~\cite{6}).
Hence, the set of $X^iY^j$, where $X=x_1 \mod \langle f_1, f_2\rangle$, 
$Y=x_2 \mod \langle f_1, f_2\rangle$, $i=1,\ldots,n-1$, $j=1,\ldots,n-2$,
is a basis of the quotient algebra.

To calculate the degree of the discriminant by Kleiman's formula,
we have to calculate\break 
$\int_Zc(X^{n-1}Y^{n-2})$. 
Let $\tilde w_0$ be the longest element in $W$, and let
$w_0$ be the shortest element in $\tilde w_0W_P$ 
with the reduced factorization
$$w_0=\left(\matrix
1&  2&3&4&\ldots&n\cr
n&n-1&1&2&\ldots&n-2\cr
\endmatrix\right)=
(n-1,n)(n-2,n-1)\ldots(12)\cdot(n-1,n)(n-2,n-1)\ldots(23).$$

Let 
$A_{w_0}=A_{(n-1,n)}A_{(n-2,n-1)}\ldots A_{(12)}
A_{(n-1,n)}A_{(n-2,n-1)}\ldots A_{(23)}$ be the corresponding endomorphism of degree
$-(2n-3)$ in $S$, where
$A_{(ij)}={\displaystyle id-s_{(ij)} \over\displaystyle x_i-x_j}$, and
$s_{(ij)}$ is the reflection that transposes
$x_i$ and $x_j$.
Then
$\int_Zc(X^{n-1}Y^{n-2})=A_{w_0}(x_1^{n-1}x_2^{n-2})$ (see~\cite{1}).
It is obvious that
$$A_{w_0}(x_1^{n-1}x_2^{n-2})=
A_{\rho_1}A_{\rho_2}(x_1^{n-1}x_2^{n-2})=
A_{\rho_1}(x_1^{n-1}A_{\rho_2}(x_2^{n-2}))=
A_{\rho_1}(x_1)A_{\rho_2}(x_2),$$
where $\rho_k=(n-1,n)(n-2,n-1)\ldots(k,k+1)$.
These factors are both equal to~$1$, since they are equal to
$\int_{\Bbb P^{n-1}}c_1({\cal O}(1))^{n-1}$ and
$\int_{\Bbb P^{n-2}}c_1({\cal O}(1))^{n-2}$, respectively,
and it is obvious that these integrals are equal to~1.
We finally obtain that $\int_Zc(X^{n-1}Y^{n-2})=1$.

It remains to calculate the polynomial
$$\sum_{i=0}^{2n-3}(i+1)c_{2n-3-i}(aX+Y)^i,\eqno(10)$$
in the ring
$\Bbb Q[X,Y]$ (with the basis $X^iY^j$, $i=1,\ldots,n-1$, $j=1,\ldots,n-2$,
and relations $X^{n-1}+X^{n-2}Y+\ldots+Y^{n-1}=0$, $Y^n=0$, and
$X^n=0$, which follows from the preceding relations), where
$c_k$ is the $k$th homogeneous component of the polynomial
$$(1-X+Y)\prod_{i=3}^n(1-X+x_i)
\prod_{i=3}^n(1-Y+x_i),$$
in which the $i$th symmetric function of
$x_3,\ldots,x_n$ must be replaced by 
$(-1)^i(X^i+X^{i-1}Y+\ldots+Y^i)$.
The result of this calculation is
$\deg(\D)X^{n-1}Y^{n-2}$.

Note that the polynomial (10) is equal to
$$F'(T)|_{T=2X+Y}=(F_1F_2F_3F_4)'(T)|_{T=2X+Y},$$
where
$$\aligned
F_1=T,\quad F_2=T-X+Y,\cr
F_3=\prod_{i=3}^n(T-X+x_i)=\sum_{i=0}^{n-2}(T-X)^{n-2-i}(-1)^i(X^i+\ldots+Y^i),\cr
F_4=\prod_{i=3}^n(T-Y+x_i)=\sum_{i=0}^{n-2}(T-Y)^{n-2-i}(-1)^i(X^i+\ldots+Y^i).
\endaligned$$
We have, further,
$$\aligned
F_3(T)=
\left(\sum_{i=0}^{n-2}(T-X)^{n-2-i}(-1)^i
{X^{i+1}-Y^{i+1}\over X-Y}\right)=\cr
{X(T-X)^{n-2}\over X-Y}\left(\sum_{i=0}^{n-2}(T-X)^{-i}(-1)^iX^{i}\right)
-{Y(T-X)^{n-2}\over X-Y}
\left(\sum_{i=0}^{n-2}(T-X)^{-i}(-1)^iY^{i}\right)=\cr
{X\left((T-X)^{n-1}-(-X)^{n-1}\right)\over T(X-Y)}-
{Y\left((T-X)^{n-1}-(-Y)^{n-1}\right)\over (T-X+Y)(X-Y)}.
\endaligned$$
We obtain, likewise, that
$$\aligned
F_4(T)=
\left(\sum_{i=0}^{n-2}(T-Y)^{n-2-i}(-1)^i
{X^{i+1}-Y^{i+1}\over X-Y}\right)=\cr
{X(T-Y)^{n-2}\over X-Y}
\left(\sum_{i=0}^{n-2}(T-Y)^{-i}(-1)^iX^{i}\right)
-{Y(T-Y)^{n-2}\over X-Y}
\left(\sum_{i=0}^{n-2}(T-Y)^{-i}(-1)^iY^{i}\right)=\cr
{X\left((T-Y)^{n-1}-(-X)^{n-1}\right)\over (T+X-Y)(X-Y)}-
{Y\left((T-Y)^{n-1}-(-Y)^{n-1}\right)\over T(X-Y)}.
\endaligned$$
We deduce from the latest formula that
$$\aligned 
F'(T)=
{n(X-Y)(T-X)^{n-1}+(-X)^n-(-Y)^n\over X-Y}\times
\qquad\qquad\qquad\qquad\qquad\qquad\cr
{(X-Y)(T-Y)^n+T(-X)^n-(T+X-Y)(-Y)^n\over T(T+X-Y)(X-Y)}+\cr
{(X-Y)(T-X)^n+(T-X+Y)(-X)^n-T(-Y)^n\over T(X-Y)}\times
\qquad\qquad\qquad\qquad\qquad\qquad\cr
{\biggl(n(T+X-Y)-(T-Y)\biggr)(T-Y)^{n-1}+(-X)^n 
\over (T+X-Y)^2}-\cr
{(X-Y)(T-X)^n+(T-X+Y)(-X)^n-T(-Y)^n\over T(X-Y)}\times
\qquad\qquad\qquad\qquad\qquad\qquad\cr
{(X-Y)(T-Y)^n+T(-X)^n-(T+X-Y)(-Y)^n\over T(T+X-Y)(X-Y)}.
\endaligned$$
Using the elementary formulae
$$\aligned
{\alpha^n(\gamma-\beta)+\beta^n(\alpha-\gamma)+\gamma^n(\beta-\alpha)\over
(\alpha-\beta)(\beta-\gamma)(\gamma-\alpha)}
=\sum_{i+j+k=n-2}\alpha^i\beta^j\gamma^k,\cr
{\alpha^n(\gamma-\beta)+\beta^n(\alpha-\gamma)+\gamma^n(\beta-\alpha)\over
(\alpha-\beta)(\gamma-\beta)}
=\sum_{i=0}^{n-2}\beta^i(\alpha^{n-1-i}-\gamma^{n-1-i}),\cr
{(n(\alpha-\beta)-\alpha)\alpha^{n-1}+\beta^{n}\over
(\alpha-\beta)^2}=\sum_{i=0}^{n-2}(n-1-i)\alpha^{n-2-i}\beta^i,
\endaligned$$
we obtain
$$\aligned
F'(T)=\biggl(n(T-X)^{n-1}+(-1)^n\sum_{i=0}^{n-1}X^{n-1-i}Y^i\biggr)\times
\biggl(\sum_{i+j+k=n-2}(T-Y)^i(-X)^j(-Y)^k\biggr)+\cr
\biggl(\sum_{i=0}^{n-2}(-X)^i((T-X)^{n-1-i}-(-Y)^{n-1-i})\biggr)\times
\biggl(\sum_{i=0}^{n-2}(n-1-i)(T-Y)^{n-2-i}(-X)^i\biggr)-\cr
\biggl(\sum_{i=0}^{n-2}(-X)^i((T-X)^{n-1-i}-(-Y)^{n-1-i})\biggr)\times
\biggl(\sum_{i+j+k=n-2}(T-Y)^i(-X)^j(-Y)^k\biggr).
\endaligned$$
Putting $T=aX+bY$ in the latest formula, we obtain
$$\aligned
\biggl(\sum_{i=0}^{n-1}
\bigl(1+(-1)^nn(a-1)^{n-1-i}b^i{n-1\choose i}\bigr)X^{n-1-i}Y^i\biggr)\times
\biggl(\sum_{i=0}^{n-2}\schoose{n-2}{i}_{a,b-1}X^iY^{n-2-i}\biggr)-\cr
\biggl(\sum_{i=0}^{n-1}\bigl(\schoose{n-1}{i}_{a-1, b}-
\schoose{n-2}{i-1}_{a-1, b}-1\bigr)X^{n-1-i}Y^i\biggr)\times
\qquad\qquad\qquad\qquad\cr
\biggl(\sum_{i=0}^{n-2}
\bigl(\fchoose{n-2}{i}_{a, b-1}-\fchoose{n-3}{i}_{a, b-1}\bigr)
X^{n-2-i}Y^i\biggr)+\cr
\biggl(\sum_{i=0}^{n-1}\bigl(\schoose{n-1}{i}_{a-1, b}-
\schoose{n-2}{i-1}_{a-1, b}-1\bigr)X^iY^{n-1-i}\biggr)\times
\biggl(\sum_{i=0}^{n-2}\schoose{n-2}{i}_{a, b-1}X^iY^{n-2-i}\biggr),
\endaligned$$
where
$${\schoose nk}_{xy}=
\sum_{p=0}^k\sum_{q=0}^{n-k}(-1)^{p+q}x^py^q{p+q\choose p},\qquad
{\fchoose nk}_{xy}=
\sum_{p=0}^k\sum_{q=0}^{n-k}(p+q+1)(-1)^{p+q}x^py^q{p+q\choose p}.
$$

It remains to calculate the degree of the discriminant, that is,
the difference between the coefficients of
$X^{n-1}Y^{n-2}$ and $X^{n-2}Y^{n-1}$ in the expression
for $F'(aX+bY)$. After some transformations we obtain
$$\aligned
\sum_{i=0}^{n-1}
\bigl((-1)^nn(a-1)^{n-1-i}b^i{n-1\choose i}
+\schoose{n-1}{n-1-i}_{a-1, b}-\schoose{n-2}{n-1-i}_{a-1, b}\bigr)\times
\qquad\qquad\qquad\cr
\bigl(\schoose{n-2}{i}_{a,b-1}-\schoose{n-2}{i-1}_{a,b-1}\bigr)-\cr
\sum_{i=0}^{n-1}
\bigl(\schoose{n-1}{n-1-i}_{a-1, b}-\schoose{n-2}{n-1-i}_{a-1, b}-1\bigr)
\times\qquad\qquad\qquad\qquad\qquad\qquad
\qquad\qquad\qquad\cr
\bigl(\fchoose{n-2}{i}_{a, b-1}-\fchoose{n-3}{i}_{a, b-1}
-\fchoose{n-2}{i-1}_{a, b-1}+\fchoose{n-3}{i-1}_{a, b-1}\bigr).
\endaligned$$
This is the degree of the discriminant of the irreducible
$\SL_n$-module with the highest weight $(a-b)\varphi_1+b\varphi_2$.

Substituting $b=1$ in the last formula, we obtain
$${(n^2-n)a^{n+1}-(n^2+n)a^{n-1}-2n(-1)^n\over (a+1)^2}.$$
If $a=2$, then we get
$${(3n^2-5n)2^n-4n(-1)^n\over 18},$$
which completes the proof of Theorems 2 and 2'.\qed

\head \S3. $E$-regular algebras\endhead

\subhead Proof of Theorem 3\endsubhead 
The proof of assertion (ii) is similar to that of assertion (iv)
of Theorem~1.

We claim that assertion (i) follows from assertion (iii).
We choose in $V^*$ a basis $\{f_1,\ldots,f_n\}$
dual to the basis $\{e_1,\ldots,e_n\}$ in $V$.
Let $U\in\Gr(n-1,V)$ be the hyperplane $f_1=0$.
It is clear that $\tilde\E=\GL_n\cdot M_0$, where
$M_0=\tilde\E\cap (U, P\A)$.
Moreover, $M=Z\cap (U, P\A)$ is the linear subspace of the algebras
for which
$U$ is an $(n-1)$-dimensional subalgebra.
Let $P$ be the parabolic subgroup
of matrices
$\left(\matrix A&0\cr *&B\cr\endmatrix\right)$,
and let $\goth u$ be the Lie algebra of matrices
$\left(\matrix0&X\cr 0&0\cr\endmatrix\right)$, where
$B$ is an $(n-1)\times(n-1)$-matrix and
$X$ is an $1\times(n-1)$-matrix.
Every algebra $A\in M$ defines alinear map
$\goth u\to\Lambda^{n-1}U^*\otimes V/U$.
Since $A\in M_0$ if and only if this map is degenerate,
we have $\codim_MM_0=1$.
It is obvious that $M_0$ is irreducible, since
$M_0$ is the spreading of the subspace
$$M_0^1=\{A\in M_0\,|\,\text{the algebra}\ E_{12}A
\ \text{has a subalgebra}\  U\}$$
by the group $P$.
Hence, $\tilde\E$ is an irreducible divisor in 
$Z$, and assertion (i) follows from assertion (iii).

We prove (iii).
We shall say that an $(n-1)$-dimensional subalgebra $U'$
of $A\in\E$ is {\it critical\/} if $(U', A)$ lies in $\tilde\E$.
In this case there is an $(n-2)$-dimensional
subspace $W'\subset U'$ such that if 
$V=U'\oplus\Bbb Ce$ and $v\in\gl_n$ is a non-zero linear operator
such that  $v(V)\subset \Bbb Ce$ and $v(W')=0$, 
then $U'$ is a subalgebra of $vA$.
(See the proof of assertion (v) of Theorem~1.)
To prove assertion (iii) it is sufficient to prove
that in generic algebras in
$M_0^1$ the subalgebra $U$ is the unique
critical subalgebra. Let
$N\subset M_0^1$ be the subvariety of all algebras that have another 
critical subalgebra.

Note that $M_0^1$ is normalized by the parabolic subgroup (9).
We denote it by $Q$. Then $N$ is the spreading of the subvarieties
$N_2$ and $N_3$ by the group $Q$, where $A\in N_i$ 
if and only if the hyperplane $f_i=0$
is a critical subalgebra.
In turn, $N_2$ is the spreading of the vector spaces
$N_2^1$ and $N_2^3$, and
$N_3$ is the spreading of the subspaces $N_3^1$, $N_3^2$, and $N_3^4$, 
where $N_i^j\subset N_i$ is the subspace of all algebras
such that the $(n-2)$-dimensional subspace $W'$ (mentioned above)
is given by $f_i=f_j=0$. 
Let $Q_i^j\subset Q$ be the subgroup that normalizes the flag
$f_i\subset\langle f_i,f_j\rangle$.
It is easy to verify that
$$\codim_QQ_2^1=1,\ \codim_QQ_2^3=n-1,\ \codim_QQ_3^1=n-1,
\ \codim_QQ_3^2=n,\ \codim_QQ_3^4=2n-3.$$

On the other hand,
$$\codim_{M_0^1}N_2^1=n,\ 
\codim_{M_0^1}N_2^3=\codim_{M_0^1}N_3^1=\codim_{M_0^1}N_3^2
=\codim_{M_0^1}N_3^4=2n-2.$$
Hence, $\codim_{M_0^1}QN_i^j\ge1$ in all cases,
which completes the proof of the theorem.\qed

\head \S4. Regular $4$-dimensional anticommutative algebras\endhead
\subhead Proof of Theorem 4 \endsubhead                              
We begin with assertion (i).
We have toprove that if $A$ is a $4$-dimensional
regular anticommutative algebra with zero trace,
then the set of its three-dimensional subalgebras $S_1,S_2,\ldots,S_5$ 
is a generic configuration of hyperplanes, that is,
the intersection of any three of them
is one-dimensional
and the intersection
of any four of them
is zero-dimensional.
Indeed, assume, for example,
that $U=S_1\cap S_2\cap S_3$ is two-dimensional.
Let $v\in U$ and $v\ne0$. Then $[v,\cdot]$ induces a linear operator on
$A/U$, since $U$ is a subalgebra.
$S_1/U$, $S_2/U$,and $S_3/U$ are one-dimensional eigenspaces.
Since $\dim A/U=2$, the operator is a homothety.
Since this is true for any $v\in U$, any three-dimensional
subspace that contains $U$ is a three-dimensionalsubalgebra,
which contradicts the fact
that there are precisely five such subalgebras.

Now assume that $U=S_1\cap S_2\cap S_3\cap S_4$ is
one-dimensional, and let $v\in U$, $v\ne0$.
Then the operator $[v,\cdot]$ induces an operator on $A/U$. 
This operator has four two-dimensional eigenspaces
$S_1/U,\ldots,S_4/U$ of which any three have zero intersection.
Hence,this operator is a homothety.
Let $W\supset U$ be an arbitrary two-dimensional subspace,
and let  $w\in W$ be a vector that is not proportional to $v$. 
Then the operator $[w,\cdot]$ induces alinear operator on $A/W$. 
Let $z$ be a non-zero eigenvector.
Then  $\langle v,w,z\rangle$ is a three-dimensional subalgebra.
Therefore, every vector can be included in a three-dimensional subalgebra,
which contradicts the fact
that there are only five such subalgebras.

This argument also shows that $A$ has no one-dimensional ideals.
Since every three-dimensional subspace that contains a two-dimensional ideal
is a subalgebra, there are no two-dimensional ideals, which completes the proof
of assertion (ii).

To prove assertion (iii), we consider the subvariety
$X$ of two-dimensional subalgebras in $A$. 
Then $X\subset\Gr(2,4)$.
Consider the Pl\"ucker embedding 
$\Gr(2,4)\subset\Bbb P^5=P(\Lambda^2\Bbb C^4)$.
First we claim that the embedding $X\subset\Bbb P^5$ 
is non-degenerate, that is, the image is contained in no hyperplane.
Let $S_1,S_2,S_3,S_4$ be four three-dimensional subalgebras.
Since it is a generic configuration, we can choose
a basis
$\{e_1,e_2,e_3,e_4\}$ in $A$ such that
$S_i=\langle e_1,\ldots,\hat e_i,\ldots,e_4\rangle$. 
Since the intersection of three-dimensional subalgebras
is a two-dimensional subalgebra,
$A$ has six subalgebras
$\langle e_i, e_j\rangle$, $i\ne j$. 
The set of corresponding bivectors
$e_i\wedge e_j$ is a basis in $\Lambda^2\Bbb C^4$.
Therefore, they can lie in no hyperplane.

Hence, $H^0(X, {\cal O}_X(1))=\Lambda^2\Bbb C^4$. 
To prove that  $X$ is a del Pezzo surface
of degree five, we have only to verify
that ${\cal O}_X(1)$ coincides with the anticanonical sheaf (see~\cite{3}).
Since $Y=\Gr(2,4)$ is a quadric in $\Bbb P^5$,
we have $\omega_Y={\cal O}_Y(2-5-1)={\cal O}_Y(-4)$.
The set $X$~is a non-singular subvariety of codimension~$2$
in~$Y$. Therefore, $\omega_X=\omega_Y\otimes\Lambda^2{\cal N}_{X/Y}$,
where ${\cal N}_{X/Y}$ is the normal sheaf.
Further, $X$ is the scheme of zeros
of a regular section of the fibre bundle
${\cal L}=\Lambda^2\SS^*\otimes V/\SS$
(see the proof of Theorem~1), whence
${\cal N}_{X/Y}={\cal L}|_Y$ and
$\Lambda^2{\cal N}_{X/Y}={\cal O}_X(3)$, since
$c_1({\cal L})=3H$. We obtain that
$\omega_X={\cal O}_X(-4)\otimes{\cal O}_X(3)={\cal O}_X(-1)$, 
as was to be shown.

It remains to prove assertion (iv).
Since $X$ is a del Pezzo surface of degree five,
it contains ten straight lines.
Since the embedding $X\subset P(\Lambda^2\Bbb C^4)$ is anticanonical,
these straight lines
are ordinary straight lines in
$P(\Lambda^2\Bbb C^4)$ that lie in $\Gr(2,4)$. 
It remains to establish a bijection between these straight lines
and fans.
If $b\in\Lambda^2\Bbb C^4$, then
$b$ belongs to the cone over $\Gr(2,4)$ if and only if
$b\wedge b=0$. If $b_1$ and $b_2$ belong to this cone,
then the straight line that joins them
belongs to this cone
if and only if $b_1\wedge b_2=0$,
which coincides with the fan condition.\qed

\widestnumber\key{XX}
\head Bibliography\endhead

\ref\key 1
\by I.N. Bernstein, I.M. Gel'fand, S.I. Gel'fand
\paper Schubert cells and cohomology of the spaces $G/P$
\jour Uspechi Mat. Nauk
\vol 28/3
\yr 1973
\pages 3--26. English transl., Russian Math. Surveys, 28/3, 1973, 1--26
\endref

\ref\key 2
\by S.L. Kleiman
\paper The enumerative theory of singularities
\jour Uspechi Mat. Nauk
\vol 35/6
\yr 1980
\pages 69--148. English transl., Russian Math. Surveys, 
35/6, 1980, 69--148
\endref

\ref\key 3
\by Yu.I. Manin
\book Cubic forms. Algebra, geometry, arithmetic
\publ M., Nauka
\yr 1972. English transl., North--Holland, Amsterdam--London 1974
\endref

\ref\key 4
\by D. Mumford
\book Complex projective varieties
\publ Springer--Verlag, Berlin--Heidelberg--New York
\yr 1979
\endref

\ref\key 5
\by W. Fulton
\book Intersection theory
\publ Springer--Verlag, Berlin--Heidelberg--New York
\yr 1984
\endref

\ref\key 6
\by G.M. Bergman
\paper The diamond lemma for ring theory
\jour Adv. in Math.
\vol 29/2
\yr 1978
\pages 178--218
\endref

\ref\key 7
\by R. Bott
\paper Homogeneous vector bundles
\jour Ann. Math.
\vol 66/2
\yr 1957
\pages 203--248
\endref

\ref\key 8
\by C. De Concini, J. Weyman
\paper A formula with nonnegative terms for the degree
of the dual variety of a homogeneous space
\jour Proceedings of the AMS
\vol 125/1
\yr 1997
\pages 1--8
\endref

\ref\key 9
\by F. Knop, G. Menzel
\paper Duale Variet\"aten von Fahnenvariet\"aten
\jour Comment. Math. Helv.
\vol 62/1
\yr 1987
\pages 38--61
\endref

\ref\key 10
\by A. Lascoux
\paper Degree of the dual Grassman variety
\jour Comm. Algebra
\vol 9/11
\yr 1981
\pages 1215--1225
\endref

\ref\key 11
\by E.A. Tevelev
\paper Generic Algebras
\jour Transformation Groups
\vol 1
\moreref Nos. 1\&2
\yr 1996
\pages 127--151
\endref

\enddocument